\documentclass[11pt]{amsart}
\usepackage{pxfonts}
\usepackage{amsmath}
\usepackage{amssymb}
\usepackage{amsthm}
\usepackage{amscd, amsfonts}
\usepackage[dvips]{epsfig}
\usepackage{verbatim}
\usepackage{epsf}
\usepackage[usenames]{color}
\usepackage[bookmarksnumbered,pdfpagelabels=true,plainpages=false,colorlinks=true,
            linkcolor=black,citecolor=black,urlcolor=blue]{hyperref}

\newcommand{\al}{\alpha}

\newcommand{\ep}{\epsilon}

\newcommand{\pbp}{\partial \bar{\partial}}

\newcommand{\hth}{\hat{\theta}}

\theoremstyle{plain}
\newtheorem{theorem}{Theorem}[section]
\newtheorem{lemma}[theorem]{Lemma}
\newtheorem{prop}[theorem]{Proposition}
\newtheorem{coro}[theorem]{Corollary}

\theoremstyle{definition}
\newtheorem{rema}[theorem]{Remark}

\numberwithin{equation}{section}

\title[Deformed Hermitian Yang-Mills PDE]{A note on the deformed Hermitian Yang-Mills PDE}

\author[Pingali]{Vamsi P. Pingali}
\address{Department of Mathematics\\
Indian Institute of Science,\\ Bangalore, India -560012}
\email{vamsipingali@iisc.ac.in}

\begin{document}
\maketitle

\begin{abstract}
We prove \emph{a priori} estimates for a generalised Monge-Amp\`ere PDE with ``non-constant coefficients" thus improving a result of Sun in the K\"ahler case. We apply this result to the deformed Hermitian Yang-Mills (dHYM) equation of Jacob-Yau to obtain an existence result and \emph{a priori} estimates for some ranges of the phase angle assuming the existence of a subsolution. We then generalise a theorem of Collins-Sz\`ekelyhidi on toric varieties and use it to address a conjecture of Collins-Jacob-Yau. 
\end{abstract}

\section{Introduction}
\indent According to some versions of superstring theory the spacetime of the universe is constrained to be a product (more generally a fibration) of a compact Calabi-Yau three-fold and a four dimensional Lorentzian manifold. A ``duality" relates the geometry of this Calabi-Yau manifold with another ``mirror" Calabi-Yau manifold. From a differential geometry standpoint this might be thought of (roughly) as a relationship between the existence of ``nice" metrics on a line bundle on one Calabi-Yau manifold  and special Lagrangian submanifolds of the other Calabi-Yau manifold. Using the Fourier-Mukai transform on a torus fibration, Leung-Yau-Zaslow suggested \cite{leung} that this ``mirror symmetry" implies that the following special Lagrangian-type PDE ought to be satisfied on a compact K\"ahler manifold $(X,\omega)$ of complex dimension $n$. 
\begin{gather}\label{speclag}
Im \left ( (\omega - F)^n \right ) = \tan (\hat{\theta}) Re \left ( (\omega - F) ^n\right ),
\end{gather}
where $F$ is the curvature of the Chern connection of a hermitian metric on a holomorphic line bundle $L$ and $\tan (\hat{\theta})$ is the constant defined by integrating equation \ref{speclag} on both sides. The angle $\hth$ has been named ``the phase angle". Throughout this paper, if $\alpha, \beta$ are forms, then we denote $\alpha \wedge \beta$ by $\alpha \beta$. Also, $\alpha^k = \alpha \wedge \alpha \ldots$ (k times).\\
\indent In \cite{special}, Jacob and Yau showed that for a given ample line bundle $L$ over a compact K\"ahler manifold with non-negative orthogonal bisectional curvature, $L^k$ admits a solution to the equation \ref{speclag}. In \cite{Col0} Collins, Jacob and Yau showed existence under the assumption of a subsolution with supercritical phase. In this paper we aim to get similar results using a different method. \\
\indent We first prove \emph{a priori} estimates on a K\"ahler manifold for a generalised Monge-Amp\`ere PDE with ``non-constant" coefficients. This is an improvement of a result of Sun \cite{Weisun}.
\begin{theorem}\label{lapgen}
Assume that $n\Omega ^{n-1} - \displaystyle \sum _{k=0}^{n-1} {n \choose k} c_k (n-k) \Omega^{n-k-1} \omega^k > 0 $. Let $\phi$ satisfying $\inf_X \phi = 1$ be a smooth solution of
\begin{gather}
\Omega_{\phi}^n = \displaystyle \sum _{k=0}^{n-1} {n \choose k} c_k \Omega_{\phi}^{n-k} \omega^k, \label{geneq}
\end{gather}
where $c_k \geq 0$ are smooth functions such that either $c_k=0$ or $c_k >0$ throughout, and $\sum _k c_k >0$. Then $\Vert \phi \Vert _{C^{2,\alpha}} \leq C (c_k, \Omega, \omega, X, \alpha)$.
\end{theorem}
\indent We then apply the preceding result to equation \ref{speclag} to arrive at the following result.
\begin{theorem}\label{mainthmhigh}
Let $(L,h_0)$ be a hermitian holomorphic line bundle on a compact K\"ahler $n$-complex dimensional manifold $(X,\omega)$. Let $F_0$ denote the curvature of the Chern connection of $(L,h_0)$. Let $\hth$ be a smooth function on $X$. Under the following assumptions any smooth solution $h=h_0e^{-\phi}$ of equation \ref{speclag} satisfies the \emph{a priori} estimate $\Vert \phi \Vert _{C^{2,\alpha}} \leq C(\alpha, X, \Omega, \hat{\theta}, \omega)$. In addition, if $\hth$ is a constant such that the integrals on both sides of equation \ref{speclag} are equal, then there exists a smooth solution satisfying \ref{speclag}.
\begin{enumerate}
\item If $n=2m+1$ then $\tan(\hth)$ is assumed to be smooth and $\Omega = \sqrt{-1}F_0-\omega \tan(\hth)$ is assumed to be a supercritical phase subsolution, i.e., $\Omega >0$ and $n\Omega^{n-1}-\displaystyle \sum c_k k\Omega ^{k-1} \omega^{n-k}>0$ where $c_k$ are defined as 
    \begin{gather}
    c_{2j} =(-1)^{m+1+j} {2m+1 \choose 2j} \sec(\hth)^{2m+2-2j}\sin((2m-2j)\hth)  \nonumber\\
c_{2j+1} =  (-1)^{m+j+1} {2m+1 \choose 2j+1} \sec(\hth)^{2m-2j+1}\cos((2m-2j-1)\hth).\label{oddexp}
    \end{gather}
    In addition, we assume that $c_k$ as defined above are either strictly positive throughout or are identically zero.
\item If $n=2m$ then $\cot(\hth)$ is assumed to be smooth and $\Omega = \sqrt{-1}F_0+\omega \cot(\hth)$ is assumed to be a supercritical phase subsolution, i.e., $\Omega >0$ and $n\Omega^{n-1}-\displaystyle \sum c_k k\Omega ^{k-1} \omega^{n-k}>0$ where $c_k$ are defined as
    \begin{gather}
    c_k = \csc(\hth) ^{2m-k} (-1)^{2m-k+1} {2m \choose k} \sin((2m-k-1)\hth).\label{evenexp}
    \end{gather}
    In addition, we assume that $c_k$ as defined above are either strictly positive throughout or are identically zero.
\end{enumerate}
\end{theorem}
\begin{rema}\label{crit}
The assumption $\Omega>0$ in theorem \ref{mainthmhigh} is implied by the assumption of supercritical phase ($\sum \arctan (\mu_l) = \hat{\theta} \in ((n-2)\frac{\pi}{2}, n\frac{\pi}{2})$ where $\mu_l$ are the eigenvalues of $\sqrt{-1}F$ with respect to $\omega$) using proposition 8.1 and lemma 8.2 of \cite{Col0}.	 
\end{rema}
\begin{rema}\label{expla}
Theorem \ref{mainthmhigh} is actually a special case of theorems 1.1 and 1.2 of \cite{Col0}. But the point here is not so much the theorem itself as the proof of it. The estimates involved and the method of continuity used are much simpler compared to those in \cite{Col0}. 
\end{rema}
\indent The main assumption of having a subsolution in theorem \ref{mainthmhigh} is actually quite difficult to check. It seems that it is almost as difficult as actually solving the equation itself. Therefore, it is desirable to have easier ``algebro-geometric" assumptions on a generalised Monge-Amp\`ere equation akin to the J-flow equation in \cite{Coll}. Indeed, by modifying the technique in \cite{Coll} we get the following result on toric varieties. Notice that the following result is also an improvement of \cite{Coll} even in spirit because we are allowing $c_1=0$ in this theorem (at the cost of assuming that $c_n>0$).
\begin{theorem}\label{torusthm}
Let $X$ be a compact toric manifold of dimension $n$, with two toric K\"ahler metrics $\Omega$ and $\omega$. Suppose that
\begin{gather}
\displaystyle \int _X \Omega^n - \sum _{k=1}^{n} \int_X {n \choose k} c_k \omega ^k \Omega ^{n-k} \geq 0, \label{first}
\end{gather}
where $c_k \geq 0$ are constants such that either $c_1>0$ or $c_n>0$. Also assume that
\begin{gather}
 \int _V \Omega ^{p} -\sum _k\int_V  c_k {p\choose k}  \omega ^k \Omega ^{p-k} > 0
\label{integ}
\end{gather}
for all $p$-dimensional toric subvarieties and all $p \leq n-1$. Then there is a smooth K\"ahler metric $\Omega _{\phi} = \Omega + \sqrt{-1}\partial \bar{\partial} \phi$ such that
\begin{gather}
\Omega _{\phi}^n = \displaystyle \sum _{k=1}^{n} c_k {n \choose k} \omega ^k \Omega_{\phi} ^{n-k} + d \omega ^n \label{colleq}
\end{gather}
for a suitable constant $d\geq0$ with equality holding if and only if equality holds in inequality \ref{first}.
\end{theorem}
Using this theorem we get the following corollary which addresses conjecture  1.4 of \cite{Col0}. 
\begin{coro}\label{corotorusthm}
For every toric subvariety $V$ of a toric manifold $X$, define $\theta _V = \mathrm{Arg} \displaystyle \int _V (\omega -F)^{dim \ V}$. Assume that $\omega$ and $F_0$ are torus invariant. Let $\hat{\theta}$ as in equation \ref{speclag} be a constant. If
\begin{enumerate}
\item For every toric subvariety of $X$, $\theta _V > \hat{\theta} - (n-\mathrm{dim} V)\frac{\pi}{2}$, and
\item $c_k$ as defined in theorem \ref{mainthmhigh} satisfy $c_k \geq 0$ with $c_0>0$,
\end{enumerate}
then there exists a smooth solution to deformed Hermitian Yang-Mills equation \ref{speclag}.
\end{coro}
\indent Almost all the results in this paper are based on the observation that equation \ref{speclag} can be treated as a generalised Monge-Amp\`ere equation. We believe that one can push these techniques further in order to obtain more general results. This will be investigated in future work.\\

\emph{Acknowledgements} : The author acknowledges the support of an SERB grant : ECR/2016/001356. He also thanks the Infosys foundation for the Infosys young investigator award. The author is deeply indebted to the anonymous referee for a very thorough reading of the paper and for useful comments.

\section{Some \emph{a priori} estimates for the generalised Monge-Amp\`ere equation}
In this section we prove a fairly general \emph{a priori} estimate on compact K\"ahler n-dimensional manifolds $(X,\omega)$ with another K\"ahler metric $\Omega _{\phi}=\Omega + \sqrt{-1} \pbp \phi$. First one has a uniform estimate.
\begin{lemma}\label{uniform}
Assume that $n\Omega ^{n-1} - \displaystyle \sum _{k=0}^{n-1} {n \choose k} c_k (n-k) \Omega^{n-k-1} \omega^k > 0 $ where $\Omega>0$ is a K\"ahler metric. Let $\phi$ satisfying $\sup_X \phi = 0$ be a smooth solution of
\begin{gather}
\Omega_{\phi}^n = \displaystyle \sum _{k=0}^{n-1} {n \choose k} c_k \Omega_{\phi}^{n-k} \omega^k,
\end{gather}
where $c_k \geq 0$ are smooth functions and $\Omega_{\phi}>0$. Then $\Vert \phi \Vert _{C^0} \leq C (c_k, \Omega, \omega, X)$.
\end{lemma}
The proof of this lemma follows almost word-to-word from \cite{Sz} and is hence omitted. Alternatively one can use the ``standard" Moser iteration type argument (\cite{Weisun} for instance). \\
\indent Next we prove a Laplacian bound for generalised Monge-Amp\`ere equations with ``non-constant coefficients" (an improvement of a result of Sun (see \cite{Weisun}) in the K\"ahler case). We change the normalisation of $\phi$ to be $\inf \phi =1$. For this estimate we need the uniform estimate.\\

\emph{Proof of theorem \ref{lapgen}} :  Define $\sigma_k(A)$ as the coefficient of $t^k$ in $\det(I+tA)$ and $S_k(A)=\frac{\sigma_{n-k}(A)}{\sigma_n(A)}$ for Hermitian matrices $A$ whose eigenvalues will be denoted by $\lambda_i$. Recall that for diagonal positive matrices $A$ and arbitrary Hermitian matrices $B$ we have the following result (lemma 8 in \cite{Coll}). We denote all constants by $C$ in whatever follows.
\begin{lemma}
Consider the function $S_k(A)$ on the cone of positive-definite Hermitian matrices $A$. Then $A \rightarrow S_k(A)$ is convex. In fact, if $B_{i\bar{j}}$ is any Hermitian matrix and $A$ is diagonal with eigenvalues $\lambda_i >0$, then we have
\begin{gather}
\displaystyle \sum \frac{\partial^2 S_k}{\partial a_{i\bar{j}} \partial a_{k\bar{l}}}(A) B_{i\bar{j}} B_{\bar{k}l} + \sum \frac{\partial S_k}{\partial  \lambda_i}(A)\frac{\vert B_{i\bar{j}} \vert^2}{\lambda _j} \geq 0. \label{propsS}
\end{gather}
\label{standardresult}
\end{lemma}
From now onwards, $S_k(\Omega_{\phi})$ will be denoted simply as $S_k$ for the sake of brevity. Akin to Yau's celebrated proof of the Calabi conjecture, let $\psi = e^{-\gamma(\phi)}\frac{n\omega^{n-1}\Omega_{\phi}}{\omega^n}$ be a function where $\gamma (x)$ is to be chosen later. Assume that the maximum of $\psi$ is attained at $p$ and that coordinates are chosen near $p$ so that $\omega (p) = \sqrt{-1} \sum dz^i \wedge d\bar{z}^i$, $\omega _{,i} (p) = \omega_{,\bar{i}}(p) = 0$, and $\Omega_{\phi} = \sqrt{-1} \sum \lambda _i dz^i \wedge d\bar{z}^i$. Here $\omega_{,\bar{i}}=\frac{\partial \omega}{\partial \bar{z}^i}$ and likewise for $\omega_{,i}$. From now onwards, a comma used as a subscript indicates a partial derivative with respect to the chosen coordinate (\emph{not} a covariant derivative). At $p$, $\nabla \psi = 0$. Using this we have,
\begin{gather}
0=\psi_{,i}(p)=\frac{n\omega^{n-1}\Omega_{\phi,i}}{\omega^n}-\gamma ^{'} \phi_{,i} \frac{n\omega^{n-1}\Omega_{\phi}}{\omega^n} 
\label{firstderivativezero}
\end{gather}
In whatever follows, we heavily exploit the fact that $\frac{\partial S_k}{\partial \lambda_i}<0$. Moreover, $\nabla ^2 \psi(p)$ is negative-definite. Hence,
\begin{gather}
0\geq \displaystyle \sum _{k,i}(-c_k \frac{\partial S_k}{\partial \lambda _i}) \psi_{,i\bar{i}}(p) \nonumber \\
\Rightarrow 0 \geq \displaystyle \sum _{k,i}(-c_k \frac{\partial S_k}{\partial \lambda _i}) \Big [\frac{n\omega^{n-1}\Omega_{\phi,i}}{\omega^n}-\gamma ^{'} \phi_{,i} \frac{n\omega^{n-1}\Omega_{\phi}}{\omega^n} \Big ] (-\gamma ^{'})\phi_{,\bar{i}} \nonumber \\
+  \displaystyle \sum _{k,i}(-c_k \frac{\partial S_k}{\partial \lambda _i}) \Big [\frac{n(n-1)\omega^{n-2}\omega_{,i\bar{i}} \Omega_{\phi}}{\omega^n} - \frac{n\omega^{n-1}\Omega_{\phi}}{\omega^n}\frac{n\omega^{n-1}\omega_{,i\bar{i}}}{\omega^n} + \frac{n\omega^{n-1}\Omega_{\phi,i\bar{i}}}{\omega^n}\nonumber\\-\gamma^{''} \vert \phi _{,i} \vert^2 \frac{n\omega^{n-1}\Omega_{\phi}}{\omega^n}-\gamma^{'} \phi_{,i\bar{i}} \frac{n\omega^{n-1}\Omega_{\phi}}{\omega^n}-\gamma^{'} \phi_{,i} \frac{n\omega^{n-1}\Omega_{\phi,\bar{i}}}{\omega^n}\Big ]. \label{diffpsi}
\end{gather}
At this point we divide equation \ref{geneq} on both sides by $\Omega_{\phi}^n$ and differentiate once with respect to $z^{\mu}$ to obtain at $p$
\begin{gather}
0 = \displaystyle \sum _{k} c_{k,\mu} S_k + \sum _{k,i} c_k \frac{\partial S_k}{\partial \lambda_i} (\Omega_{\phi})_{i\bar{i},\mu} \ \forall  \ \mu.\label{diffonce}
\end{gather}
Differentiating again (with respect to $\bar{z}^{\mu}$), summing over $\mu$, and using inequality \ref{propsS} we get
\begin{gather}
0\geq \displaystyle \sum_{k,\mu} c_{k,\mu\bar{\mu}} S_k + \sum_{k,\mu} c_{k,\mu} \frac{\partial S_k}{\partial \lambda _i} (\Omega_{\phi})_{i\bar{i},\bar{\mu}} + \sum_{k,\mu} c_{k,\bar{\mu}} \frac{\partial S_k}{\partial \lambda _i} (\Omega_{\phi})_{i\bar{i},\mu} + \sum_{k,\mu} c_k {n \choose k} \frac{\Omega_{\phi} ^{n-k} k \omega^{k-1} \omega_{,\mu \bar{\mu}}}{\Omega_{\phi}^n} \nonumber \\
+ \sum_{k,\mu,i} c_{k} \frac{\partial S_k}{\partial \lambda _i} (\Omega_{\phi})_{i\bar{i},\mu\bar{\mu}} - \sum_{k,\mu,i,j} c_{k} \frac{\partial S_k}{\partial \lambda _i} \frac{\vert(\Omega_{\phi})_{i\bar{j},\mu} \vert^2}{\lambda_j}.\label{difftwice}
\end{gather}
Note that $\frac{n\omega^{n-1}\Omega_{\phi,i\bar{i}}}{\omega^n}(p)=\displaystyle \sum_{\mu} (\Omega_{\phi})_{i\bar{i},\mu\bar{\mu}}$. With this in mind, using \ref{difftwice} in \ref{diffpsi} we get 
\begin{gather}
0 \geq \displaystyle \sum_{k,\mu} c_{k,\mu\bar{\mu}} S_k + \sum_{k,\mu} c_{k,\mu} \frac{\partial S_k}{\partial \lambda _i} (\Omega_{\phi})_{i\bar{i},\bar{\mu}} + \sum_{k,\mu} c_{k,\bar{\mu}} \frac{\partial S_k}{\partial \lambda _i} (\Omega_{\phi})_{i\bar{i},\mu} + \sum_{k,\mu} c_k {n \choose k} \frac{\Omega_{\phi} ^{n-k} k \omega^{k-1} \omega_{,\mu \bar{\mu}}}{\Omega_{\phi}^n} \nonumber \\
- \sum_{k,\mu,i,j} c_{k} \frac{\partial S_k}{\partial \lambda _i} \frac{\vert(\Omega_{\phi})_{i\bar{j},\mu} \vert^2}{\lambda_j} +  \displaystyle \sum _{k,i}(-c_k \frac{\partial S_k}{\partial \lambda _i}) \Big [\frac{n\omega^{n-1}\Omega_{\phi,i}}{\omega^n}-\gamma ^{'} \phi_{,i} \frac{n\omega^{n-1}\Omega_{\phi}}{\omega^n} \Big ] (-\gamma ^{'})\phi_{,\bar{i}} \nonumber \\
+  \displaystyle \sum _{k,i}(-c_k \frac{\partial S_k}{\partial \lambda _i}) \Bigg [\frac{n(n-1)\omega^{n-2}\omega_{,i\bar{i}} \Omega_{\phi}}{\omega^n} - \frac{n\omega^{n-1}\Omega_{\phi}}{\omega^n}\frac{n\omega^{n-1}\omega_{,i\bar{i}}}{\omega^n} \nonumber\\-\gamma^{''} \vert \phi _{,i} \vert^2 \frac{n\omega^{n-1}\Omega_{\phi}}{\omega^n}-\gamma^{'} \phi_{,i\bar{i}} \frac{n\omega^{n-1}\Omega_{\phi}}{\omega^n}-\gamma^{'} \phi_{,i} \frac{n\omega^{n-1}\Omega_{\phi,\bar{i}}}{\omega^n}\Bigg ] \nonumber 
\end{gather}
At this point we note that $\vert c_{k,\mu \bar{\mu}} \vert \leq C$ and $0<\displaystyle \sum_{k : c_k>0} S_k \leq C$ (here we are using the assumption that either $c_k>0$ everywhere or $=0$ identically). Also, $-C\omega \leq \omega_{\mu \bar{\mu}} \leq C\omega$, $-C\omega \leq \omega_{\mu} \leq C\omega$, and $-C\omega \leq \omega_{\bar{\mu}} \leq C\omega$.  Therefore, upon grouping terms and simplifying we get the following.
\begin{gather}
C\geq \displaystyle \sum _{k,i}(-c_k \frac{\partial S_k}{\partial \lambda _i}) \Big [\frac{n\omega^{n-1}\Omega_{\phi,i}}{\omega^n} (-\gamma ^{'} \phi_{,\bar{i}})+\frac{n\omega^{n-1}\Omega_{\phi,\bar{i}}}{\omega^n} (-\gamma ^{'} \phi_{,i})+((\gamma^{'})^2-\gamma^{''}) \vert \phi _{,i} \vert^2 \frac{n\omega^{n-1}\Omega_{\phi}}{\omega^n}\nonumber \\-\gamma^{'}\phi_{,i\bar{i}}\frac{n\omega^{n-1}\Omega_{\phi}}{\omega^n} - C \frac{n\omega^{n-1}\Omega_{\phi}}{\omega^n} \Big ] \nonumber \\
+  \sum c_{k,\mu} \frac{\partial S_k}{\partial \lambda _i} (\Omega_{\phi})_{i\bar{i},\bar{\mu}} + \sum c_{k,\bar{\mu}} \frac{\partial S_k}{\partial \lambda _i} (\Omega_{\phi})_{i\bar{i},\mu} - \sum_{k} c_{k} \frac{\partial S_k}{\partial \lambda _i} \frac{\vert(\Omega_{\phi})_{i\bar{j},\mu} \vert^2}{\lambda_j}. \label{subs}
\end{gather}
At this point we notice that (here  we again use the assumption that either $c_k$ is identically zero or $c_k >0$ everywhere)
\begin{gather}
\sum c_{k,\bar{\mu}} \frac{\partial S_k}{\partial \lambda _i} (\Omega_{\phi})_{i\bar{i},\mu} +  \sum c_{k,\mu} \frac{\partial S_k}{\partial \lambda _i} (\Omega_{\phi})_{i\bar{i},\bar{\mu}} \geq \displaystyle \sum c_k \frac{\partial S_k}{\partial \lambda _i} \Bigg ( \epsilon \frac{\vert (\Omega_{\phi})_{i\bar{i},\mu}\vert^2}{\lambda_i} + \frac{1}{\epsilon} \vert(\ln  c_k)_{,\bar{\mu}}\vert^2 \lambda_i   \Bigg ) \nonumber \\
= \displaystyle \sum c_k \frac{\partial S_k}{\partial \lambda _i} \epsilon \frac{\vert (\Omega_{\phi})_{i\bar{i},\mu}\vert^2}{\lambda_i} + \sum \frac{1}{\epsilon} \vert(\ln  c_k)_{,\bar{\mu}}\vert^2 c_k k S_k \nonumber \\
\geq \displaystyle \sum c_k \frac{\partial S_k}{\partial \lambda _i} \epsilon \frac{\vert (\Omega_{\phi})_{i\bar{i},\mu}\vert^2}{\lambda_i} - \frac{C}{\epsilon} . \label{herewe}
\end{gather}
Likewise, using the Cauchy-Schwarz inequality on the first two terms of \ref{subs} we see that
\begin{gather}
\frac{n\omega^{n-1}\Omega_{\phi,i}}{\omega^n} (-\gamma ^{'} \phi_{,\bar{i}})+\frac{n\omega^{n-1}\Omega_{\phi,\bar{i}}}{\omega^n} (-\gamma ^{'} \phi_{,i}) \geq - \frac{\omega^n}{n\omega^{n-1}\Omega_{\phi}}\vert \frac{n\omega^{n-1}\Omega_{\phi,i}}{\omega^n} \vert^2 - (\gamma^{'}) ^2 \vert \phi_{,i} \vert^2   \frac{n\omega^{n-1}\Omega_{\phi}}{\omega^n}.
\label{firsttwo}
\end{gather}
Putting \ref{subs}, \ref{herewe}, and \ref{firsttwo} together we obtain
\begin{gather}
C+\frac{C}{\epsilon} \geq \sum _{k,i}(-c_k \frac{\partial S_k}{\partial \lambda _i}) \Big [- \frac{\omega^n}{n\omega^{n-1}\Omega_{\phi}}\vert \frac{n\omega^{n-1}\Omega_{\phi,i}}{\omega^n} \vert^2 -(\gamma^{''}\vert \phi _{,i} \vert^2+C) \frac{n\omega^{n-1}\Omega_{\phi}}{\omega^n}-\gamma^{'}\phi_{,i\bar{i}}\frac{n\omega^{n-1}\Omega_{\phi}}{\omega^n} \Big ] \nonumber \\
+ \sum c_k \frac{\partial S_k}{\partial \lambda _i} \epsilon \frac{\vert (\Omega_{\phi})_{i\bar{i},\mu}\vert^2}{\lambda_i} - \sum c_{k} \frac{\partial S_k}{\partial \lambda _i} \frac{\vert(\Omega_{\phi})_{i\bar{j},\mu} \vert^2}{\lambda_j}.
\label{together}
\end{gather}
It is easy to see (for instance using Lagrange multipliers for the second inequality) that
\begin{gather}
\sum c_k \frac{\partial S_k}{\partial \lambda _i} \frac{\vert (\Omega_{\phi})_{i\bar{i},\mu}\vert^2}{\lambda_i} - \sum c_{k} \frac{\partial S_k}{\partial \lambda _i} \frac{\vert(\Omega_{\phi})_{i\bar{j},\mu} \vert^2}{\lambda_j} \geq 0 \nonumber \\
- \frac{\omega^n}{n\omega^{n-1}\Omega_{\phi}}\vert \frac{n\omega^{n-1}\Omega_{\phi,i}}{\omega^n} \vert^2 + \displaystyle \sum _{j,\mu} \frac{\vert(\Omega_{\phi})_{i\bar{j},\mu} \vert^2}{\lambda_j}  \geq 0.
\label{easytosee}
\end{gather}
Using \ref{easytosee} in \ref{together} we see that
\begin{gather}
C+\frac{C}{\epsilon} \geq \sum _{k,i}(-c_k \frac{\partial S_k}{\partial \lambda _i}) \Big [- \epsilon\frac{\omega^n}{n\omega^{n-1}\Omega_{\phi}}\vert \frac{n\omega^{n-1}\Omega_{\phi,i}}{\omega^n} \vert^2 -(\gamma^{''}\vert \phi _{,i} \vert^2+C) \frac{n\omega^{n-1}\Omega_{\phi}}{\omega^n}-\gamma^{'}\phi_{,i\bar{i}}\frac{n\omega^{n-1}\Omega_{\phi}}{\omega^n} \Big ].
\label{aftereasytosee}
\end{gather}
Using equation \ref{firstderivativezero} we see that
\begin{gather}
C+\frac{C}{\epsilon} \geq \sum _{k,i}(-c_k \frac{\partial S_k}{\partial \lambda _i}) \Big [ -((\epsilon(\gamma^{'})^2+\gamma^{''})\vert \phi _{,i} \vert^2+C) \frac{n\omega^{n-1}\Omega_{\phi}}{\omega^n}-\gamma^{'}\phi_{,i\bar{i}}\frac{n\omega^{n-1}\Omega_{\phi}}{\omega^n} \Big ].
\label{aftereasytoseetwo}
\end{gather}
Now we use a standard lemma (originally due to Fang-Lai-Ma \cite{InvHess}) in the following form due to Wei Sun (lemma 4.1 in \cite{Weisun}).
\begin{lemma}
There are constants $N,\theta>0$ such that when $\Delta \phi(p)=\displaystyle \sum_i \phi_{,i\bar{i}}>N$, at $p$ we have the following inequality.
\begin{gather}
\displaystyle \sum_{k,i} -c_k \frac{\partial S_k}{\partial \lambda_i}\phi_{,i\bar{i}}\leq -\theta(\displaystyle \sum_{k,i}-c_k \frac{\partial S_k}{\partial \lambda_i} +1).
\end{gather}
\label{standardlemma}
\end{lemma}
Using lemma \ref{standardlemma} we see that (assuming that $\gamma^{'}\geq0$)
\begin{gather}
C+\frac{C}{\epsilon} \geq \sum _{k,i}(-c_k \frac{\partial S_k}{\partial \lambda _i}) \Big [ -((\epsilon(\gamma^{'})^2+\gamma^{''})\vert \phi _{,i} \vert^2+C) \frac{n\omega^{n-1}\Omega_{\phi}}{\omega^n} \Big ] + \theta \gamma^{'}\frac{n\omega^{n-1}\Omega_{\phi}}{\omega^n}  \Big [ 1+\sum -c_k \frac{\partial S_k}{\partial \lambda _i} \Big ],
\end{gather}
for some uniform constant $\theta>0$. Choosing $\gamma (x) = \frac{\ln(x)}{\epsilon}$ we see that $-\gamma^{''}-\epsilon(\gamma^{'})^2 = 0$, and $\gamma^{'} \geq \frac{1}{\epsilon \sup \phi}$. Therefore,
\begin{gather}
C+\frac{C}{\epsilon} \geq \sum _{k,i}(-c_k \frac{\partial S_k}{\partial \lambda _i}) \Big [ -C \frac{n\omega^{n-1}\Omega_{\phi}}{\omega^n} \Big ] + \frac{\theta}{\epsilon \sup \phi} \frac{n\omega^{n-1}\Omega_{\phi}}{\omega^n}  \Big [ 1+\sum -c_k \frac{\partial S_k}{\partial \lambda _i} \Big ].
\end{gather}
 Choosing $\epsilon<\frac{\theta}{C\sup\phi}$ we see that $\frac{n\omega^{n-1}\Omega_{\phi}}{\omega^n} (p)$ has to be bounded above to avoid a contradiction. Since $p$ is the point of maximum of $\psi=e^{-\gamma(\phi)}\frac{n\omega^{n-1}\Omega_{\phi}}{\omega^n}$ and $\Vert \phi \Vert_{C^0} \leq C$, this means that $\psi \leq C$ on the manifold. Thus, $\Delta \phi \leq C$ everywhere. \qed \\
\textbf{}\\

Before we prove theorem \ref{mainthmhigh}, we prove the following lemma on higher order \emph{a priori} estimates.
\begin{lemma}
If $c_k$ are constants, then $\Vert \phi \Vert_{C^{k,\alpha}} \leq C$ for all $k$.
\label{higherorderestimates}
\end{lemma}
\begin{proof}
In \cite{Siuevans} a complex version of the Evans-Krylov theory was developed. It was proven if $F(u_{i\bar{j}})=f(x)$ where $F$ is uniformly elliptic (i.e., $CId \geq \frac{\partial F}{\partial u_{i\bar{j}}} \geq C^{-1}Id$) and concave, then $\Vert \Delta u \Vert_{C^{0,\alpha}} \leq C$ where $C$ depends on the ellipticity constants of $F$ (which may in turn depend on $\Vert\Delta u \Vert_{C^0}$) and $\Vert f \Vert_{C^2}$). In our case, we rewrite equation \ref{geneq} as 
\begin{gather}
-1=\displaystyle \sum_{k=0}^{n-1} {n \choose k}\frac{-c_k \Omega_{\phi}^{n-k}\omega^k}{\Omega_{\phi}^n}
\label{rewrittengeneq}
\end{gather}
It is a well known fact that the functions $A \rightarrow -S_k(A)$ are elliptic and concave on the cone of positive-definite matrices (indeed, this is implied by lemma \ref{standardresult} for instance). The Laplacian bound we proved in theorem \ref{lapgen} is easily seen to imply uniform bounds on the ellipticity constants. Thus by the complex version of the Evans-Krylov theory we can conclude that $\Vert \Delta u \Vert_{C^{0,\alpha}}\leq C$. By the Schauder estimates, this implies that $\Vert u \Vert_{C^{2,\alpha}}\leq C$. Differentiating again and bootstrapping this regularity using the Schauder estimates we see that $\Vert \phi \Vert_{C^{k,\alpha}} \leq C$. 
\end{proof}
Armed with this lemma we proceed further. \\

\emph{Proof of theorem \ref{mainthmhigh}} : Using the method of continuity as in \cite{myap} one can (re)prove an existence result for the generalised Monge-Amp\`ere PDE in the case where $c_k$ are constants \cite{myap,Weisun,Coll}. We rewrite equation \ref{speclag} as follows.\\
 If $n=2m+1$, let $\Omega = \sqrt{-1}F_0-\omega \tan(\hth)$.
\begin{gather}
\mathrm{Im}((\omega+\sqrt{-1}(\Omega_{\phi}+\omega \tan(\hth)))^{2m+1}) = \tan(\hth) \mathrm{Re}((\omega+\sqrt{-1}(\Omega_{\phi}+\omega \tan(\hth)))^{2m+1}) \nonumber \\
\Rightarrow \displaystyle \sum _{k=0}^{2m+1} {2m+1 \choose k} \omega^{2m+1-k} \Omega_{\phi} ^{k} \mathrm{Im} ((\sqrt{-1})^k(1+\sqrt{-1}\tan(\hth))^{2m+1-k}) \nonumber \\
 =  \tan(\hth)\displaystyle \sum _{k=0}^{2m+1} {2m+1 \choose k} \omega^{2m+1-k} \Omega_{\phi} ^{k} \mathrm{Re} ((\sqrt{-1})^k(1+\sqrt{-1}\tan(\hth))^{2m+1-k}) \nonumber \\
 \Rightarrow \displaystyle \sum _{k=0}^{2m+1} {2m+1 \choose k} \sec(\hth)^{2m+1-k} \omega^{2m+1-k} \Omega_{\phi} ^{k} \mathrm{Im} ((\sqrt{-1})^k(\cos((2m+1-k)\hth)+\sqrt{-1}\sin((2m+1-k)\hth))) \nonumber \\
 =  \tan(\hth)\displaystyle \sum _{k=0}^{2m+1} {2m+1 \choose k} \sec(\hth)^{2m+1-k} \omega^{2m+1-k} \Omega_{\phi} ^{k} \mathrm{Re} ((\sqrt{-1})^k(\cos((2m+1-k)\hth)+\sqrt{-1}\sin((2m+1-k)\hth))) 
\label{speclaggenmaodd}
\end{gather} 
 If $n=2m$, let $\Omega = \sqrt{-1}F_0+\omega \cot(\hth)$. Calculations similar to \ref{speclaggenmaodd} show that
 \begin{gather}
\displaystyle \cot(\hth) \sum _{k=0}^{2m} (-1)^k{2m \choose k} \csc(\hth)^{2m-k} \omega^{2m-k} \Omega_{\phi} ^{k} \sin((2m-k)\hth) \nonumber \\
 =  \displaystyle \sum _{k=0}^{2m} (-1)^k {2m \choose k} \csc(\hth)^{2m-k} \omega^{2m-k} \Omega_{\phi} ^{k}\cos((2m-k)\hth) 
\label{speclaggenmaeven}
\end{gather} 
In either case, \ref{speclag} boils down to
\begin{gather}
\Omega_{\phi}^n = \displaystyle \sum c_k \Omega_{\phi}^{k} \omega ^{n-k},
\label{speclagma}
\end{gather}
where if $n=2m+1$
\begin{gather}
c_k =(-1)^m {2m+1 \choose k} \sec(\hth)^{2m+1-k} \Bigg ( \tan(\hth)\mathrm{Re} ((\sqrt{-1})^k(\cos((2m+1-k)\hth)+\sqrt{-1}\sin((2m+1-k)\hth))) \nonumber \\
 - \mathrm{Im} ((\sqrt{-1})^k(\cos((2m+1-k)\hth)+\sqrt{-1}\sin((2m+1-k)\hth))) \Bigg ), 
\end{gather}
and if $n=2m$
\begin{gather}
c_k = \csc(\hth) ^{2m-k} (-1)^k {2m \choose k} \Bigg [ \cot(\hth) \sin((2m-k)\hth)-\cos((2m-k)\hth) \Bigg ].
\end{gather}
\indent At this point we recall the method of continuity used in \cite{myap}. 
\begin{gather}
\Omega_{\phi_t}^n=t\displaystyle \sum_{k=1}^{n-1} c_k \Omega_{\phi_t}^k \omega^{n-k}+c^{1-t}b_{t}c_0\omega^n,
\label{mycontinuitypath}
\end{gather}
where $t\in[0,1]$, $c= \frac{\displaystyle \int \Omega^n}{\displaystyle \int c_0\omega^n}$, and $b_t$ is a normalising constant chosen so that the integrals are equal on both sides, i.e., $b_t =c^{t-1}\frac{\displaystyle \int ((1-t)\Omega^n+tc_0\omega^n)}{\displaystyle \int c_0 \omega^n}$ (hence $b_t c^{1-t} \geq 1$). At $t=0$ the equation corresponds to the Calabi conjecture solved by Yau. Openness was proven in \cite{myap}. The \emph{a priori} estimates proven earlier show closedness and hence the equation can be solved for $t=1$. Unfortunately, $c_0$ is allowed to be zero and hence the above argument can break down in that case. \\
\indent Suppose $c_0=0$. In this case,  we set up an approximate version of \ref{speclagma}.
\begin{gather}
\Omega_{\phi_{\epsilon}}^n=\displaystyle \sum_{k=1}^{n-1}c_k a_{\epsilon} \Omega_{\phi_{\epsilon}}^k \omega^{n-k}+\epsilon \omega^n,
\label{approximateversion}
\end{gather}
where $a_{\epsilon}= \frac{\displaystyle \int \Omega^n - \epsilon \int \omega^n}{\int \Omega^n}$. By the argument above, we can solve \ref{approximateversion} for all $\epsilon>0$. Note that $a_{\epsilon} \rightarrow 1$ as $\epsilon \rightarrow 0$. Therefore, since all the relevant quantities in the proof of theorem \ref{lapgen} remain bounded, the \emph{a priori} estimates proven in theorem \ref{lapgen} show that $\Vert \phi_{\epsilon}\Vert_{C^{k,\alpha}} \leq C_k$ where $C_k$ is independent of $\epsilon$. Hence, a subsequence of $\phi_{1/n}$ converges in $C^{k,\beta}$ (for $\beta<\alpha$ and $k>2$) to a function $\phi$ solving \ref{speclagma}. Since this holds true for all $k$ and uniqueness holds for $C^2$ solutions, we see that $\phi$ is smooth. 
\section{Toric varieties}
We first state and prove a theorem that is almost the same as theorem \ref{torusthm} with one extremely important difference, namely, $c_1$ is assumed to be positive. The proof of the following theorem follows \cite{Coll} very closely (we just use more general $C^{2,\al}$ estimates borrowed from \cite{pin}).
\begin{theorem}\label{pretorusthm}
Assuming the same hypotheses as in theorem \ref{torusthm} and that $c_1>0$ we have the same conclusion as in theorem \ref{torusthm}, i.e., the generalised Monge-Amp\`ere equation has a solution.
\end{theorem}
\emph{Proof of theorem \ref{pretorusthm}} : Just as in \cite{Coll} the proof proceeds by induction on the dimension $n$ of $X$. For one dimensional manifolds the result is straightforward. Consider the method of continuity
\begin{gather}
t= \displaystyle \sum _{k=1}^{n} c_k {n \choose k} \frac{\omega ^k \Omega_{\phi_t} ^{n-k}}{\Omega_{\phi_t}^n} + d_t \frac{\omega ^n}{\Omega_{\phi_t}^n}.
\label{newmethodcont}
\end{gather}
For $t\rightarrow \infty$ there exists a solution using the Calabi-Yau theorem and the implicit function theorem. Thus the problem reduces to proving that the infimum of $t\geq 0$ for which the equation has a solution is $0$. For this one needs $C^{2,\alpha}$ \emph{a priori} estimates on $\phi$. Using the induction hypothesis one sees that the equation
\begin{gather}
t\Omega _{\phi_i}^{n-1} = \displaystyle \sum _{k=1}^{n-1} c_k {n-1 \choose k} \omega ^k \Omega_{\phi_i} ^{n-1-k} + r_i \omega ^{n-1}
\label{lowerdim}
\end{gather}
has a smooth torus-invariant solution over each of the toric divisors $D_i$ of $X$ where $r_i\geq 0$ is some constant.\\
\indent At this point, using the $\phi_i$, just as in \cite{Coll} one can construct a viscosity supersolution $\phi$ to the equation in a neighbourhood $U$ of the union of the toric divisors $D=\cup_i D_i$. The definition of a viscosity supersolution to equation \ref{newmethodcont} is somewhat technical. Suppose $B$ is a positive-definite Hermitian $n\times n$ matrix. Let $\lambda_1,\ldots, \lambda_n$ be its eigenvalues and $D_n=diag(\lambda_1,\ldots,\lambda_n)$.  Now let $\tilde{F}(B) = \max \displaystyle \sum_{k=1}^{n-1} c_k S_k(D_{n-1})$ where the maximum runs over all $n-1$-tuples of eigenvalues of $B$.\\
\indent  A K\"ahler current $\chi$ with a continuous local potential on $U\subset X$ is said to be a strict viscosity supersolution of \ref{newmethodcont} if $\chi=\Omega+\sqrt{-1}\pbp \psi$ satisfies the following : If $p \in U$, and $h : V \rightarrow \mathbb{R}$ is a $C^2$ function on a neighbourhood $V$ of $p$ where $\chi = \sqrt{-1}\pbp f$, then whenever $h-f$ has a local minimum at $p$, then $\tilde{F}(\omega^{i\bar{q}}\partial_j \partial_{\bar{q}}h) \leq t-\delta$ at $p$ for some positive $\delta$ (independent of $p$).
\\

\indent The point is that such a supersolution is a ``barrier" function that allows one to prove estimates. In particular, one has the following estimate.
\begin{prop}[A special case of proposition 16 in \cite{Coll}]
Suppose $\Omega_{\phi}$ is smooth and satisfies equation \ref{newmethodcont} on $\bar{U} \subset X$. Suppose also that we have a strict viscosity supersolution $\chi =\Omega+\sqrt{-1}\pbp \psi$ of \ref{newmethodcont}. Then we have an estimate $\frac{n\omega^{n-1}\Omega_{\phi}}{\omega^n} \leq Ce^{N(\phi-\inf \phi)}$ where $C$ may depend on the maximum of $\frac{n\omega^{n-1}\Omega_{\phi}}{\omega^n}$ on the boundary of $U$.
\end{prop}
\indent The construction of such a supersolution will be done in the following lemma.
\begin{lemma}
There exists a neighbourhood $U$ of $\cup_i D_i$ and a strict viscosity supersolution $\chi=\Omega+\sqrt{-1}\pbp \psi$ on $U$.
\end{lemma}
\begin{proof}
\indent Extend $\phi_i$ to all of $X$. Then define $\chi_i = \Omega_{\phi_i}+A\sqrt{-1}\pbp (\gamma(d_i) \vert d_i \vert^2)$ where $d_i$ is the distance to $D_i$ and $\gamma : \mathbb{R}\rightarrow \mathbb{R}_{\geq 0}$ is a cutoff function supported near $0$ (and equal to $1$ on a smaller neighourhood of $0$ where the distance function is smooth). Choose $A>>1$ so that $\chi_i >0$ on a neighbourhood $U_i$ of $D_i$ and satisfies $G(\chi_i)=\displaystyle \sum _{k=1}^{n} c_k {n \choose k} \frac{\omega ^k \chi_i^{n-k}}{\chi_i^n}<t-\epsilon$ for some small $\epsilon>0$. Write $\chi_i =\Omega+\sqrt{-1}\pbp\psi_i$. Define $\tilde{\psi}_i=\psi_i -B_i +\delta \displaystyle \sum_{j<i}\gamma(d_j) \ln d_j$. Choose the constants so that $\tilde{\chi}_i=\Omega+\sqrt{-1}\pbp\tilde{\psi}_i>0$ and satisfies $G(\tilde{\chi}_i)<t-\epsilon$ on a neighbourhood of $D_i \cap \cap_{j<i}D_j^c$. Also, we choose $B_i$ inductively so that for $j>i$, on $D_j$ we have $\psi_i-B_i<\psi_j-B_j$. Now define $\psi = \max_i \tilde{\psi}_i$. Lemma 15 of \cite{Coll} implies that the maximum of supersolutions is a supersolution. (Recall that since according to our definition, $\tilde{F}(B)$ is minus of an elliptic operator in the usual sense, ``supersolution" in our sense means ``subsolution" as in the usual sense (like $\Delta u \geq 0$) for elliptic equations. Hence, the maximum of our supersolutions is a supersolution.) Using this observation it can be seen that $\psi$ is the desired supersolution after shrinking $U$.  
\end{proof}
\indent If one manages to prove a $C^{2,\alpha}$ estimate on the torus-invariant $\phi_t$ outside $U$, then one can use the fact due to Tosatti and Weinkove \cite{ToWe} that proving the following estimate $$\frac{n\omega^{n-1}\Omega_{\phi}}{\omega^n} \leq Ce^{N(\phi-\inf \phi)}$$ implies a $C^0$ estimate on $\phi$ and hence a global Laplacian bound on $\phi$.  
The complex version of the Evans-Krylov theorem then gives an estimate on $\Vert \Delta \phi \Vert_{C^{0,\alpha}}$. Then the Schauder estimates furnish $C^{2,\alpha}$ estimates thus completing the proof. So the problem has been reduced to proving $C^{2,\alpha}$ estimates for $\phi$ outside $U$.\\
\indent To this end, one just needs to prove such estimates on every compact set $K$ outside $U$. Since we are dealing with torus invariant metrics, we may assume that locally $\Omega_{\phi} = \sqrt{-1}\partial \bar{\partial} g$ and $\omega = \sqrt{-1}\partial \bar{\partial} f$ where $g,f :\mathbb{R}^n \rightarrow \mathbb{R}$ are convex. (The torus invariance forces the existence of such $g$ and $f$ depending only on $\vert z_i \vert$ and not on the phase angles of these complex numbers.) Now we have the following proposition akin to lemma 26 in \cite{Coll}. 
\begin{prop}
$\Vert g \Vert_{C^1(K)} \leq C$.
\end{prop}
\begin{proof}
Using an affine change we may assume that $g(0)=\nabla g(0)=0$ (where $\nabla$ is the Euclidean gradient in the chosen toric coordinates). Now the image of $\nabla g$ (being the moment map) is a convex polytope determined by the K\"ahler class of $\Omega$. Thus $\vert \nabla g(x) \vert \leq C \ \forall  \ x\in K$ and hence $\vert g(x) \vert \leq C \ \forall \ x \in K$.
\end{proof}
 Given a Laplacian bound on $\phi$ (which is equivalent to saying that $\omega^{n-1}\Omega_{\phi} \leq C\omega^n$), one can obtain a $C^{2,\alpha}$ estimate. Given that, the usual bootstrapping with Schauder estimates give us higher order estimates. Indeed, 
\begin{prop}
If $\omega^{n-1}\Omega_{\phi} \leq C\omega^n$ on $K$, then $\Vert \phi \Vert_{C^{2,\alpha}} \leq C$ (where $C$ is independent of $t$). 
\end{prop}
\begin{proof}
If we manage to prove that $\Vert \Delta_{\omega} \phi \Vert_{C^{0,\alpha}}\leq C$, then the Schauder estimate gives us a $C^{2,\alpha}$ bound. As before, the complex version of Evans-Krylov theory can be made to apply if we prove that the operator $F(A) = \displaystyle \sum_k -c_k S_k(A)$ is concave and uniformly elliptic. We already know that $F(A)$ is concave on the cone of positive-definite matrices. Hence, uniform ellipticity is all that matters. \\
\indent Since $\Omega_{\phi}>0$, and upper bound on the Laplacian on $\phi$ implies an upper bound on $\Omega_{\phi}$. Since $c_1>0$, from equation \ref{newmethodcont} we see that since $t\leq T$ for some $T$, $\Omega_{\phi} \geq C^{-1}\omega$. Thus $C^{-1} \omega \leq \Omega_{\phi} \leq C\omega$. This is easily seen to imply uniform ellipticity of $F$. This completes the proof.   
\end{proof}
\indent To prove a Laplacian bound, closely following proposition 27 of \cite{Coll}, one uses a contradiction argument in conjunction with the Legendre transform. Very roughly speaking, the Legendre transform converts $D^2 g$ to $D^2 h =(D^2 g)^{-1}$. Therefore, an upper bound on $D^2 g$ translates into a lower bound on $D^2 h$. If this is violated and $D^2 h$ becomes degenerate somewhere, then a constant-rank theorem of Bian-Guan \cite{BG} ensures that it is degenerate everywhere (this is one of the places where we use the hypothesis that $c_1>1$) and this provides the desired contradiction. Indeed,
\begin{lemma}
Suppose $f,g : B\rightarrow \mathbb{R}$ are smooth convex functions on the unit ball, $\inf_B g = g(0)= 0$ satisfying
\begin{gather}
t=\displaystyle \sum _{k=1}^{n} c_k {n \choose k} \frac{(\sum_{i,j}f_{ij}dx^i \wedge dx^j) ^k (\sum_{i,j}g_{ij}dx^i \wedge dx^j) ^{n-k}}{(\sum_{i,j}g_{ij}dx^i \wedge dx^j)^n} + d_t \frac{\det(f_{ij})}{\det(g_{ij})}.
\label{eqtoric}
\end{gather} 
Then there is a $C$ depending on $\Vert g \Vert_{C^0}$, bounds on $c_k, d_t, t, \Vert f \vert_{C^{3,\alpha}}$, and lower bounds on the Hessian of $f$ such that $\sup_{B(1/2)} \vert D^2g \vert <C$.
\end{lemma}
\begin{proof}
Assume $c=1$ by scaling $f$. As in \cite{Coll}, we use a contradiction argument. Suppose we have sequences $f_k, g_k$, satisfying $\vert g_k \vert <N$ and the other hypotheses, but $\vert D^2 g_k (x_k) \vert>k$ for some $x_k \in B(1/2)$. Then we have a uniform lower bound on the Hessians of the $g_k$ from equation \ref{eqtoric}. (Here we are using that $c_1>0$.)\\
\indent Let $h_k$ be the Legendre transform of $g_k$. Differentiating equation \ref{eqtoric} once, observing that the equation is uniformly elliptic, and using the Krylov-Safonov theorem (the nondivergence form of the De Giorgi-Nash-Moser theorem), we see that on a slightly smaller ball, $\Vert g_k \Vert_{C^{1,\alpha}}\leq C$. Hence $g_k \rightarrow g$ in $C^{1,\beta}$ (upto a subsequence) where $0<\beta<\alpha$  on a smaller ball. It is easily seen that $g$ is strictly convex. Lemma 28 in \cite{Coll} implies that $\nabla g(0.9B) \subset U_k$. For the convenience of the reader, we state it here : 
\begin{lemma}[Lemma 28, \cite{Coll}]
Suppose $f_k : B\rightarrow \mathbb{R}$ are convex, with $f_{k,ij}>\tau \delta_{ij}$, such that they converge uniformly to $f:B\rightarrow \mathbb{R}$. If $B_1 \subset B_2 \subset B_3 \subset B$ are relatively compact balls, then for sufficiently large $k$ the gradient maps satisfy $\nabla f_k (B_1) \subset \nabla f(B_2) \subset \nabla f_k (B_3)$.
\label{Collemma}
\end{lemma}
 Thus $g(\nabla 0.8B)$ is at a fixed uniform distance from $\partial U_k$. Now $h_k$ satisfies $h_k(0)=\nabla h_k (0)=0$ and the equation
\begin{gather}
t=\displaystyle  \sum _{a=1}^{n-1} c_a S_{a,D^2f_{k} (\nabla h_k)} (D^2h_{k}) +d_t S_{n,D^2f_{k}(\nabla h_k)}(h_{k,ij}),
\label{newnotation} 
\end{gather}
where $S_{a,D^2f_{k}(\nabla h)} (D^2h_{k})(x)$ is the coefficient of $t^a$ in $\det(D^2 f_k (\nabla h_k (x)) +tD^2 h_k (x))$. Using lemma \ref{noevans} below along with the Schauder estimates and the Arzela-Ascoli theorem, we extract a subsequence of $f_k$ and $h_k$ converging to $f, h : \nabla g(0.8 B) \rightarrow \mathbb{R}$ satisfying
\begin{gather}
t=\displaystyle  \sum _{a=1}^{n-1} c_a S_{a,D^2f (\nabla h)} (h_{ij}) +d_t S_{n,D^2f(\nabla h)}(D^2 h)
\end{gather} 
As in \cite{Coll} there are two cases :
\begin{enumerate}
\item If $D^2 h > CId$, then $D^2 h_k >CId \ \forall \ k>>1$ on $\nabla g(0.8B)$. So $D^2 g_k <C Id$ on $x$ such that $\nabla g_k (x) \in \nabla g(0.8 B)$. But by lemma 28 in \cite{Coll}, $\nabla g_k(0.7 B) \subset \nabla g(0.8 B)$ for large $k$ and hence we have a contradiction.
\item If the $D^2 h$ is degenerate somewhere, then applying the constant rank theorem of Bian-Guan (theorem 1.1 in \cite{BG} which has been proven recently in a simpler manner by Sz\'ekelyhidi and Weinkove \cite{SW}) we see that $D^2 h$ is degenerate everywhere and hence $\displaystyle \int_{g(\nabla 0.8 B)}\det(D^2 h) =0.$ But this contradicts $\displaystyle \int_{\nabla g_k (0.7B)} \det(D^2 h_k)=Vol(0.7 B)$. It remains to be verified that indeed the hypothesis of Bian-Guan's theorem are satisfied. To check this, we state their theorem here :
\begin{theorem}[Theorem 1.1, \cite{BG}, Theorem 1.1 in \cite{SW}]
Suppose $\mathcal{S}^n$ is the space of real symmetric $n\times n$ matrices and $U$ is a domain in $\mathbb{R}^n$. Assume that $F=F(r,p,u,x) \in C^{2,1}(\mathcal{S}^n \times \mathbb{R}^n \times \mathbb{R} \times U)$, $\frac{\partial F}{\partial r_{\alpha \beta}}>0, \ \forall  \ x\in U$, and \\
\indent $F(A^{-1},p,u,x)$ is locally convex in $(A,u,x)$ for each fixed $p$ on the cone of positive-definite $A$.  $(*)$\\
If $u\in C^{2,1}(U)$ is a convex solution of $F(D^2 u, Du, u, x)=0$, then the rank of $D^2 u(x)$ is a constant in $U$. 
\end{theorem}
In our case, $F(r,p,u,x) = \displaystyle  \sum _{a=1}^{n-1} c_a S_{a,D^2f(p)} (r) +d_t S_{n,D^2f(p)}(r)+t$. The convexity of $F(A^{-1},p,u,x)$ is a standard fact mentioned earlier. We use $c_1>0$ to show that ellipticity is satisfied by $F$. Indeed, when $r_{\alpha \beta}$ is a positive-definite matrix, $\frac{\partial F}{\partial r_{\alpha \beta}} \geq c_1 f_{\alpha \beta} >0$.
\end{enumerate}
This completes the proof.
\end{proof}
\indent To apply the Bian-Guan theorem to $h$ in the above proof, one needs the following lemma (essentially a $C^{2,\alpha}$ estimate) :
\begin{lemma}\label{noevans}
Suppose that $h:B\rightarrow \mathbb{R}$ is a smooth convex function on the unit ball in $\mathbb{R}^n$ satisfying
\begin{gather}\label{legeq}
\displaystyle  \sum _{k=2}^n c_k S_{k,b_k (\nabla h)} (h_{ij}) +\sum_{i,j}a_{ij}(\nabla h)h_{ij} = 1,
\end{gather}
where $b_k$ are positive-definite matrix-valued functions, $c_i \geq 0$ are constants,  $0<\lambda <a_{ij}<\Lambda$, $a_{ij}, b_k \in C^{1,\alpha}$, and $S_{k,b_k(\nabla h)} (h_{ij})$ is the coefficient of $t^k$ in $\det(b_{k} (\nabla h (x)) +tD^2 h (x))$. Then we have $\Vert h \Vert _{C^{2,\alpha}(\frac{1}{2}B)}<C$ for a constant $C$ depending on $\Vert h \Vert _{C^1}, \lambda, \Lambda$ and $C^{1,\alpha}$ bounds for $a_{ij}$ and $b_k$.
\end{lemma}
The proof of this lemma follows from the results in \cite{pin} (which are a generalisation of those in \cite{colcomplain}) and repeating the blowup argument of proposition 29 in \cite{Coll}. We give a sketch of the proof here.
\begin{proof}
If we prove that $N_h = \sup_{x\in B} d_x \vert D^3 h(x) \vert \leq C$ we will be done. (Indeed, this would imply that $D^2 h$ cannot vary much. So if $D^2 u$ is large somewhere, then it is uniformly large everywhere. Therefore, $Dh$ has to become large which is a contradiction.) Assume that $N_h >1$. If the maximum is achieved at $x_0 \in B$, define a function $\tilde{h} = d_{x_0}^{-2}N_h ^2 h(x_0+d_{x_0}N_h^{-1}z)-A-A_i z_i$ where $A, A_i$ are chosen so that $\tilde{h}(0)=\nabla \tilde{h}(0)=0$. It is not hard to see that $\Vert \tilde{h} \Vert_{C^3(B(2^{-1}N_h))}<C$ and $\Vert \nabla \tilde{h} \Vert_{C^0} <C$ on compact subsets of $B_{N_h/4}$. \\
\indent Suppose we have a sequence of functions $h_k$ on $B$ satisfying \ref{legeq} but with $N_{h_k}>4k$. Then the rescaled functions $\tilde{h}_k$ defined on $B_{4k}(0)$ converge (after taking a subsequence) to $\tilde{h} : \mathbb{R}^n \rightarrow \mathbb{R}$ in $C^{3,\alpha}$. The limit can be easily seen to satisfy $\vert D^3 \tilde{h}(0)\vert =1$ and 
\begin{gather}
\displaystyle  \sum _{k=2}^n c_k S_{k,\tilde{b}_k} (D^2 \tilde{h}) +\sum_{i,j}\tilde{a}_{ij}\tilde{h}_{ij} = 1,
\end{gather}
where $\tilde{b}_k, \tilde{a}_{ij}$ are constants. It is also clear that $\tilde{h}$ is smooth. At this point, a Liouville theorem can be used to conclude that $\tilde{h}$ is a quadratic polynomial thus contradicting the fact that $\vert D^3 \tilde{h}(0) \vert =1$. \\
\indent The Liouville theorem alluded to is a generalisation of lemma 31 in \cite{Coll}. 
\begin{lemma}
Suppose that $u: \mathbb{R}^n \rightarrow \mathbb{R}$ is a smooth convex function satisfying 
\begin{gather}
\displaystyle  \sum _{k=2}^n c_k S_{k,b_k} (D^2 u) +\sum_{i,j}a_{ij}u_{ij} = 1,
\label{liou}
\end{gather}
where $b_k, a$ are positive-definite constant matrices. Then $u$ is a quadratic polynomial.
\end{lemma}
This lemma's proof is omitted because it is exactly the same as that of lemma 31 in \cite{Coll}. The only difference is that we need interior $C^{2,\alpha}$ estimates for solutions of \ref{liou} which are provided by theorem 1.2 in \cite{pin}. 
\end{proof}
This completes the proof of theorem \ref{pretorusthm}.
 \qed \\

Now we use a small trick to complete the proof of theorem \ref{torusthm}. \\
\textbf{}\\
\emph{Proof of theorem \ref{torusthm}} : Assume that $c_1=0$. Since there are only finitely many toric subvarieties, the assumptions of theorem \ref{torusthm} imply that for sufficiently small $\epsilon>0$ 
$$ \displaystyle \int _X \Omega^n - \int _X n\epsilon \omega \Omega ^{n-1} - \sum _{k=2}^{n-1} \int_X {n \choose k} c_k \omega ^k \Omega ^{n-k} -(c_n-C\epsilon) \int_X \omega ^n > 0$$ for some constant $C>0$, 
and likewise for all $p$-dimensional subvarieties. Therefore, by theorem \ref{pretorusthm} we see that there is a smooth function $u_{\epsilon}$ and a constant $d_{\epsilon} >0$ such that
\begin{gather}
\Omega _{u_{\ep}}^n = \displaystyle  n\epsilon \omega \Omega_{u_{\ep}} ^{n-1} + \sum _{k=2}^{n-1} c_k {n \choose k} \omega ^k \Omega_{u_{\ep}} ^{n-k} + d_{\ep} \omega ^n. 
\end{gather}
This implies that $\Omega_{u_{\ep}}$ satisfies 
\begin{gather}
n\Omega _{u_{\ep}}^{n-1} - \sum _{k=2}^{n-1} c_k {n \choose k}(n-k) \omega ^k \Omega_{u_{\ep}} ^{n-k-1} >0 \label{satcone}
\end{gather}
Taking $\Omega _{u_{\ep}}$ as the background metric instead of $\Omega$, using the \emph{a priori} estimates in lemma \ref{lapgen} (or alternatively in \cite{Coll,Weisun}), and the method of continuity as in \cite{myap} we can solve the equation. This gives the desired result.
\textbf{}\\

We are ready to prove corollary \ref{corotorusthm}. \\

\emph{Proof of corollary \ref{corotorusthm}} : In the proof of theorem \ref{mainthmhigh} we reduced the dHYM equation to a generalised Monge-Amp\`ere equation (equation \ref{speclagma}). We may deduce the desired result by applying theorem \ref{torusthm} provided the relevant hypotheses are satisfied. Lemma $8.2$ of \cite{Col0} implies that the ``integration over subvarieties" hypotheses \ref{integ} is equivalent to assumption (1) of corollary \ref{corotorusthm}. Indeed, we state lemma $8.2$ of \cite{Col0} here from which the assertion can be seen through an easy calculation. 
\begin{lemma}[Lemma 8.2, \cite{Col0}]
A $(1,1)$-form $\chi \in [\omega]$ is a subsolution to the deformed Hermitian-Yang-Mills equation if and only if, for any $1\leq p \leq n-1$, and any non-zero, simple, positive $(n-p,n-p)$ form $\beta$, we have 
\begin{gather}
Arg\left ( \frac{(\alpha+\sqrt{-1}\chi)^p\wedge \beta}{\alpha^n}\right )> \hat{\Theta}-(n-p)\frac{\pi}{2}.
\end{gather}
\end{lemma}
This completes the proof. \qed

\end{document}